\documentstyle[12pt]{amsart}
\newtheorem{theorem}{Theorem}[section]
\newtheorem{corollary}[theorem]{Corollary}
\newtheorem{remark}[theorem]{Remark}
\newtheorem{lemma}[theorem]{Lemma}
\newtheorem{proposition}[theorem]{Proposition}

\numberwithin{equation}{section}

\begin{document}

\title[Symplectic critical surfaces in K\"ahler surfaces]
{Symplectic critical surfaces in K\"ahler surfaces}

\author{Xiaoli Han, Jiayu Li}

\thanks{The research was partially supported by  NSFC }

\address{Xiaoli Han, Math. Group, The abdus salam ICTP\\ Trieste 34100, Italy.}
\email{xhan@@ictp.it}

\address{Jiayu Li, Math. Group, The abdus salam ICTP\\ Trieste 34100,
   Italy\\
   and Academy of Mathematics and Systems Sciences\\ Chinese Academy of
Sciences\\ Beijing 100080, P. R. of China. } \email{jyli@@ictp.it}


\begin{abstract}
Let $M$ be a K\"ahler surface and $\Sigma$ be a closed symplectic
surface which is smoothly immersed in $M$. Let $\alpha$ be the
K\"ahler angle of $\Sigma$ in $M$. We first deduce the
Euler-Lagrange equation of the functional
$L=\int_{\Sigma}\frac{1}{\cos\alpha}d\mu$ in the class of
symplectic surfaces. It is $\cos^3\alpha
H=(J(J\nabla\cos\alpha)^\top)^\bot$, where $H$ is the mean
curvature vector of $\Sigma$ in $M$, $J$ is the complex structure
compatible with the K\"ahler form $\omega$ in $M$, which is an
elliptic equation. We then study the properties of the equation.
\end{abstract}
\maketitle

\section{Introduction}

Suppose that $M$ is a K\"ahler surface. Let $\omega$ be the
K\"ahler form on $M$ and let $J$ be a complex structure compatible
with $\omega$. The Riemannian metric $\langle,\rangle$ on $M$ is
defined by
$$
\langle U,V \rangle =\omega(U,JV).
$$
For a compact oriented real surface $\Sigma$ which is smoothly
immersed in $M$, one defines, following \cite {CW}, the K\"ahler
angle $\alpha$ of $\Sigma$ in $M$ by
\begin{equation}\label{e1}\omega|_\Sigma=\cos\alpha d\mu_\Sigma\end{equation} where $d\mu_\Sigma$
is the area element of $\Sigma$ of the induced metric from
$\langle,\rangle$. We say that $\Sigma$ is a holomorphic curve if
$\cos\alpha \equiv 1$, $\Sigma$ is a Lagrangian surface if
$\cos\alpha \equiv 0$ and $\Sigma$ is a symplectic surface if
$\cos\alpha > 0$.

It was conjectured by Tian \cite {T} that every embedded
orientable closed symplectic surface in a compact
K\"ahler-Einstein surface is isotopic to a symplectic minimal
surface in a suitable sense.

If the K\"ahler-Einstein surface is of nonnegative scalar
curvature, a symplectic minimal surface is holomorphic. However,
if the scalar curvature is negative, there are symplectic minimal
surfaces which are not holomorphic (\cite{A}). A symplectic
minimal surface is a critical point of the area of surfaces, which
is symplectic. It may be more natural to consider directly the
critical point of the functional
$$L=\int_{\Sigma}\frac{1}{\cos\alpha}d\mu_{\Sigma},$$
in the class of symplectic surfaces in a K\"ahler surface. It is
clear that holomorphic curves minimize the functional. The
critical point of the functional is called {\it \bf a symplectic
critical surface}.

In the paper, we first calculate the Euler-Lagrange equation of
the functional $L$.

\begin{theorem}\label{crie} Let $M$ be a K\"ahler surface.
The Euler-Lagrange equation of the functional $L$ is
$$\cos^3\alpha H- (J(J\nabla\cos\alpha)^\top)^\bot=0, $$
where $H$ is the mean curvature vector of $\Sigma$ in $M$, and
$()^\top$ means tangential components of $()$, $()^\bot$ means the
normal components of $()$. We call the critical point {\bf a
symplectic critical surface}.
\end{theorem}

We checked that it is an elliptic equation. We derive an equation
for the K\"ahler angle of a symplectic critical surface in a
K\"ahler-Einstein surface.
\begin{theorem}\label{Lapan}
If $M$ is a K\"ahler-Einstein surface and $\Sigma$ is a symplectic
critical surface, then we have
\begin{eqnarray*}
\Delta\cos\alpha &=&\frac{3\sin^2\alpha -2}{\cos\alpha}|\nabla
\alpha|^2-R\cos^3\alpha\sin^2\alpha ,
\end{eqnarray*} where $R$ is the scalar curvature of $M$.
\end{theorem}

As a corollary, we see that, if the scalar curvature $R$ of the
K\"ahler-Einstein surface $M$ is nonnegative, then a symplectic
critical surface in $M$ is holomorphic.

It is not difficult to see that, a non holomorphic symplectic
critical surface in a K\"ahler surface has at most finite complex
points. Moreover, we can show that,

\begin{theorem}\label{toplocal}
Suppose that $\Sigma$ is a non holomorphic symplectic critical
surface in a K\"ahler surface $M$. Then
$$
\chi (\Sigma)+\chi
(\nu)=-P-\frac{1}{2\pi}\int_{\Sigma}\frac{|\nabla\alpha|^2}{\cos^2\alpha}d\mu,
$$
and
$$
c_1(M)([\Sigma])=-P-\frac{1}{2\pi}\int_{\Sigma}\frac{|\nabla\alpha|^2}{\cos^3\alpha}d\mu,
$$
where $\chi(\Sigma)$ is the Euler characteristic of $\Sigma$,
$\chi(\nu)$ is the Euler characteristic of the normal bundle of
$\Sigma$ in $M$, $c_1(M)$ is the first Chern class of $M$,
$[\Sigma]\in H_2(M,{\bf Z})$ is the homology class of $\Sigma$ in
$M$, and $P$ is the number of complex tangent points.

\end{theorem}

From Theorem \ref{toplocal}, we derive some global topology
properties of a symplectic critical surface.

\begin{theorem}\label{topglobal01}
Suppose that $\Sigma$ is a symplectic critical surface in a
K\"ahler surface $M$. Then
$$
\chi (\Sigma)+\chi (\nu)\geq c_1(M)([\Sigma]),
$$
the equality holds if and only if $\Sigma$ is a holomorphic curve.

\end{theorem}

Let $g$ be the genus of $\Sigma$, $I_\Sigma$ the self-intersection
number of $\Sigma$, $D_\Sigma$ the number of double points of
$\Sigma$. Then
$$
\chi(\Sigma)=2-2g,
$$
$$
\chi(\nu)=I_\Sigma-2D_\Sigma.
$$

Set
$$
C_1(\Sigma)=C_1(M)([\Sigma]),
$$
we have
\begin{theorem}\label{topglobal02}
Suppose that $\Sigma$ is a symplectic critical surface in a
K\"ahler surface $M$. Then
$$
2-2g-
c_1(\Sigma)+I_\Sigma-2D_\Sigma=\frac{1}{2\pi}
\int_{\Sigma}\frac{(1-\cos\alpha)|\nabla\alpha|^2}{\cos^3\alpha}d\mu.
$$

\end{theorem}

In forthcoming papers, we will use variational approach and flow
method to study the existence of a symplectic critical surface in
a K\"ahler surface.

It is natural to conjecture that, {\it in each homotopy class of a
symplectic surface in a K\"ahler surface, there is a symplectic
critical surface.}

As a start point of studying the gradient flow of the function
$L$, we derive the evolution equation of $\cos\alpha$ along the
flow, which implies that the symplectic property is preserved.

\section{The Euler-Lagrange equation}

Let $\{\phi_t\}_{0\leq t\leq 1}$ be a 1-parameter family of
immersions $\Sigma\rightarrow M$ such that $\phi_0=F$ and
$\Sigma_t=\phi_t(\Sigma)$ are symplectic. Also, let $X$ denote the
initial velocity vector for $\phi_t$, i.e, $X=\left.
\frac{\partial\phi_t}{\partial t}\right|_{t=0}$. We denote by
$\bar\nabla$ the covariant derivative and by $K$ the Riemannian
curvature tensor on $M$. Furthermore, $\nabla, R$ denote the
covariant derivative and the Riemannian curvature tensor of the
induced metric $g$ on the surface $\Sigma_0$.

We start from computing the first variation of the area for this
one-parameter family of surfaces, which is in fact well-known.

\begin{proposition}\label{pro1}
The variation of the area of $\Sigma_t$ is
$$\frac{\partial}{\partial t}d\mu_t|_{t=0}=(divX^\top-X\cdot H)d\mu_t$$
\end{proposition}
{\it Proof.} The proof is routine (c.f. \cite{Si}). Fix one point
$p\in \Sigma$. Let $\{x^i\}$ be a normal coordinate system for
$\Sigma$ at $p$. Around $p$ we choose the local orthonormal frame
$\{\tilde e_1, \tilde e_2, \tilde v_3, \tilde v_4\}$ on $M$ along
$\Sigma_t$ such that $\{\tilde e_1=\partial \phi_t/\partial x^i,
\tilde e_2=\partial \phi_t/\partial x^2\}$, $\{\tilde v_3, \tilde
v_4\}$ are in the tangent bundle and in the normal bundle of
$\Sigma_t$ respectively. For simplicity, we denote $\tilde
e_i(0)=\frac{\partial F}{\partial x^i}$ by $e_i$ and identify it
with $\partial_i=\frac{\partial}{\partial x^i}$, $i=1, 2$. We also
denote $\tilde v_\alpha(0)$ by $v_\alpha$, $\alpha=3, 4$.
Furthermore, we assume that $\nabla_{e_i}e_j=0$ at $p$. Suppose
that, in this frame $X$ takes the form $X=X^i e_i+X^\alpha
v_\alpha$ and $(g_t)_{ij}=\langle\frac{\partial\phi_t}{\partial
x^i}, \frac{\partial\phi_t}{\partial x^j}\rangle$. Then,
\begin{eqnarray*}
\frac{\partial g_{ij}}{\partial t}|_{t=0}&=&
\frac{\partial}{\partial t}|_{t=0}\langle \partial\phi_t/\partial
x^i, \partial\phi_t/\partial x^j\rangle=\langle\bar\nabla_{e_i} X,
e_j\rangle+\langle e_i, \bar\nabla_{e_j}
X\rangle\\&=&\langle\bar\nabla_{e_i}(X^k e_k+X^\alpha v_\alpha),
e_j\rangle+\langle e_i, \bar\nabla_{e_j}(X^k e_k+X^\alpha
v_\alpha) \rangle\\&=& X_{j,i}-X_\alpha
h^\alpha_{ij}+X_{i,j}-X_\alpha h^\alpha_{ij}.
\end{eqnarray*} It is easy to get that
\begin{eqnarray*}\frac{\partial}{\partial t}d\mu_t|_{t=0}
&=&\frac{1}{2}g^{ij}(X_{j,i}-X_\alpha
h^\alpha_{ij}+X_{i,j}-X_\alpha
h^\alpha_{ij})d\mu_t\\&=&(divX^\top-X\cdot H)d\mu.\end{eqnarray*}
\hfill Q.E.D.

\begin{theorem}\label{th1}
Let $M$ be a K\"ahler surface. The first variational formula of
the functional $L$ is, for any smooth vector field $X$ on
$\Sigma$,
\begin{eqnarray}
\delta_X L &=& -2\int_{\Sigma}\frac{X\cdot H}{\cos\alpha}
d\mu+2\int_{\Sigma}\frac{X\cdot
(J(J\nabla\cos\alpha)^\top))^\bot}{\cos^4\alpha}d\mu,
\end{eqnarray}
where $H$ is the mean curvature vector of $\Sigma$ in $M$, and
$()^\top$ means tangential components of $()$, $()^\bot$ means the
normal components of $()$. The Euler-Lagrange equation of the
functional $L$ is
\begin{equation}\label{cri}\cos^3\alpha H-(J(J\nabla\cos\alpha)^\top)^\bot=0.
\end{equation}
\end{theorem}
{\it Proof.} We use the frame that we choose in Proposition
\ref{pro1}. From the definition of K\"ahler angle (\ref{e1}) we
have
$$\cos\alpha_t=\frac{\omega(\partial\phi_t/\partial x^1,
\partial\phi_t/\partial x^2)}{\sqrt{\det(g_t)}},
$$ where $\det (g_t)$ is the determinant of the metric $(g_t)$.
So, the functional can be written as
$$ L_t=L(\phi_t)=\int_{\Sigma}\frac{\det(g_t)}{\omega(\partial\phi_t/\partial x^1,
\partial\phi_t/\partial x^2)} dx^1\wedge dx^2.
$$
Thus using Proposition \ref{pro1} we have, \allowdisplaybreaks
\begin{eqnarray*}
\frac{d}{dt}|_{t=0} L_t &=& \int_{\Sigma}\large (\frac{\partial_t
g_{ij}|_{t=0} g^{ij}}{\omega(e_1 , e_2
)}-\frac{\partial_t\omega(\partial\phi_t/\partial x^1,
\partial\phi_t/\partial x^2)|_{t=0}}{\omega^2(e_1 , e_2 )}\large )\det(g)dx^1\wedge dx^2 \\
&=& \int_{\Sigma}(\frac{2div X^{T}-2X\cdot
H}{\cos\alpha})\sqrt{\det(g)}dx^1\wedge
dx^2\\
&&-\int_{\Sigma}(\frac{\omega(\bar\nabla_{e_1} X, e_2)+
\omega(e_1, \bar\nabla_{e_2} X)}{\cos^2\alpha})dx^1\wedge dx^2\\
&=& I+II.
\end{eqnarray*} Since $\Sigma$ is closed, applying the Stokes
formula, we obtain
\begin{eqnarray*}
I &=&\int_{\Sigma}(\frac{2\langle X, e_1\rangle
\nabla_{e_1}\cos\alpha+2\langle X,
e_2\rangle\nabla_{e_2}\cos\alpha}{\cos^2\alpha} -\frac{2 X\cdot
H}{\cos\alpha}) \sqrt{\det(g)}dx^1\wedge dx^2.
\end{eqnarray*} The second term is
\begin{eqnarray*}
II &=&-\int_{\Sigma}\frac{\nabla{e_1}(\omega(X, e_2))-\omega(X,
\bar\nabla_{e_1}\bar\nabla_{e_2} F)}{\cos^2\alpha}dx^1\wedge dx^2\\
&&-\int_{\Sigma}\frac{\nabla_{e_2}(\omega(e_1, X))-\omega(
\bar\nabla_{e_2}\bar\nabla_{e_1} F, X)}{\cos^2\alpha} dx^1\wedge
dx^2\\ &=&-\int_{\Sigma}\frac{\nabla_{e_1}(\omega(X,
e_2))+\nabla_{e_2}(\omega(e_1, X))}{\cos^2\alpha}dx^1\wedge dx^2\\
&=&-2\int_{\Sigma}\frac{\omega(X,
e_2)\nabla_{e_1}\cos\alpha+\omega(e_1,
X)\nabla_{e_2}\cos\alpha}{\cos^3\alpha} dx^1\wedge dx^2,
\end{eqnarray*} where we have used the fact that $\omega$ is
parallel. In the following, we compute pointwisely so we assume
the frame is orthonormal.  Note that,
\begin{eqnarray*}
\omega(X^\top, e_2)&=&-\langle X^\top,Je_2\rangle\\
&=&-\langle X,e_1\rangle \langle e_1,Je_2\rangle\\
&=&\langle X,e_1\rangle\cos\alpha ,
\end{eqnarray*}
and
\begin{eqnarray*}
\omega(e_1,X^\top)&=&\langle X^\top,Je_1\rangle\\
&=&\langle X,e_2\rangle \langle e_2,Je_1\rangle\\
&=&\langle X,e_2\rangle\cos\alpha ,
\end{eqnarray*}
 we separate the second term into two
parts,
\begin{eqnarray*}
II &=&-2\int_{\Sigma}\frac{\langle X,
e_1\rangle\nabla_{e_1}\cos\alpha}{\cos^2\alpha}+\frac{\langle X,
e_2\rangle\nabla_{e_2}\cos\alpha}{\cos^2\alpha}d\mu\\
&&-2\int_{\Sigma} \frac{\omega(X^\bot,
e_2)\nabla_{e_1}\cos\alpha+\omega(e_1,
X^\bot)\nabla_{e_2}\cos\alpha}{\cos^3\alpha} d\mu.
\end{eqnarray*} Therefore, we obtain that,
\begin{eqnarray*}
\frac{d}{dt}|_{t=0} L_t&=&-2\int_{\Sigma}\frac{X\cdot
H}{\cos\alpha}d\mu-2\int_{\Sigma}\frac{\omega(X^\bot,
e_2)\nabla_{e_1}\cos\alpha}{\cos^3\alpha}
d\mu\\&&-2\int_{\Sigma}\frac{\omega(e_1,
X^\bot)\nabla_{e_2}\cos\alpha}{\cos^3\alpha} d\mu.
\end{eqnarray*}

Because
\begin{eqnarray*}
(J\nabla\cos\alpha)^\top&=&(Je_1\nabla_{e_1}\cos\alpha+Je_2\nabla_{e_2}\cos\alpha)^\top\\
&=&\langle Je_1,e_2\rangle e_2\nabla_{e_1}\cos\alpha+\langle
Je_2,e_1\rangle e_1\nabla_{e_2}\cos\alpha\\
&=& (e_2\nabla_{e_1}\cos\alpha-
e_1\nabla_{e_2}\cos\alpha)\cos\alpha,
\end{eqnarray*}
we have
\begin{eqnarray*}
\frac{d}{dt}|_{t=0} L_t&=& -2\int_{\Sigma}\frac{X\cdot
H}{\cos\alpha} d\mu\\
&&-2 \int_{\Sigma}\frac{\omega(X^\bot,
(J\nabla\cos\alpha)^\top)}{\cos^4\alpha} d\mu\\ &=&
-2\int_{\Sigma}\frac{X\cdot H}{\cos\alpha}
d\mu\\
&&+2\int_{\Sigma}\frac{X\cdot
(J(J\nabla\cos\alpha)^\top))^\bot}{\cos^4\alpha} d\mu.
\end{eqnarray*} This completes the proof of the theorem.

\hfill Q. E. D.

For later purpose, and to understand the equation, we express
$(J(J\nabla\cos\alpha)^\top))^\bot$ at a fixed point $p$ in local
frame. Let $\{e_1, e_2, v_1, v_2\}$ be a orthonormal frame around
$p \in \Sigma$ such that it is normal at $p$ and $\omega, J$ take
the forms (cf. \cite {CT}, \cite {CW}),
\begin{equation}\label{e2}
\omega=\cos\alpha u_1\wedge u_2+\cos\alpha u_3\wedge
u_4+\sin\alpha u_1\wedge u_3-\sin\alpha u_2\wedge u_4
\end{equation} where $\{u_1, u_2, u_3, u_4\}$ is the dual frame of
$\{ e_1,  e_2,  v_3,  v_4\}$, and
\begin{eqnarray}\label{e3}
J =\left ( \begin{array}{clcr} 0 &\cos\alpha &\sin\alpha &0\\
-\cos\alpha &0 &0 &-\sin\alpha\\ -\sin\alpha &0 &0 &\cos\alpha
\\ 0 &\sin\alpha &-\cos\alpha &0\end{array}\right).
\end{eqnarray}

 Then
\begin{eqnarray*}
(J(J\nabla\cos\alpha)^\top))^\bot &=&
(J(\cos\alpha\partial_1\cos\alpha
e_2-\cos\alpha\partial_2\cos\alpha e_1) )^\bot \\ &=&-\cos\alpha
\sin\alpha\partial_1\cos\alpha
v_4-\cos\alpha\sin\alpha\partial_2\cos\alpha v_3 \\
&=&\cos\alpha\sin^2\alpha\partial_1\alpha
v_4+\cos\alpha\sin^2\alpha\partial_2\alpha v_3.
\end{eqnarray*}
Furthermore,
\begin{eqnarray*}
\partial_1\cos\alpha &=&\omega(\bar\nabla_{e_1}e_1, e_2)+\omega(e_1, \bar\nabla_{e_1}e_2
) \\ &=&h^\alpha_{11}\langle Jv_\alpha, e_2\rangle
+h^\alpha_{12}\langle Je_1, v_\alpha\rangle \\ &=&
(h^4_{11}+h^3_{12})\sin\alpha.
\end{eqnarray*} Similarly, we can get that,
$$\partial_2\cos\alpha=(h^3_{22}+h^4_{12})\sin\alpha.$$
Set $V=\partial_2\alpha v_3+\partial_1\alpha v_4$. Then
\begin{equation}\label{v}
V=-(h^3_{22}+h^4_{12})v_3-(h^4_{11}+h^3_{12})v_4.
\end{equation} And consequently the
Euler-Lagrange equation of the function $L$ is
\begin{equation}\label{main}
\cos^2\alpha H-\sin^2\alpha V=0.
\end{equation}
or equivalently,
$$
H-\sin^2\alpha(H+V)=0.
$$
It is not difficult to see, roughly, the symbol of the equation is
\begin{eqnarray*}
\sigma:=\left(\begin{array}{lr} (1-\sin^2\alpha)\xi^2+\eta^2
&(\sin^2\alpha)\xi\eta\\
(\sin^2\alpha)\xi\eta &\xi^2+(1-\sin^2\alpha)\eta^2
\end{array}\right),
\end{eqnarray*}
which enables one believes that the equation (\ref{cri}) is an
elliptic equation.

In the following, we give a detailed proof.
\begin{theorem}\label{th2}
The equation (\ref{cri}) is elliptic.
\end{theorem}
{\it Proof.} Assume that $\Sigma$ is immersed in $M$ by $F$.  Let
$\{x, y\}$ be a coordinate system around $p\in\Sigma$. Since
$\Sigma$ is smooth, by implicit theorem we can write $\Sigma$ as
the graph of two functions $f, g$ in a small neighborhood $U$ of
$p$, i.e, $F=(x, y, f(x, y), g(x, y))$ in $U$. Suppose that the
complex structure in the neighborhood of $F(p)$ is standard, i.e,
\begin{eqnarray*}
J=\left(\begin{array}{lccr} 0 &-1 &0 &0\\ 1 &0 &0 &0 \\ 0 &0 &0 &-1\\
0 &0 &1 &0\end{array}\right).
\end{eqnarray*} We choose $e_1=\frac{\partial F}{\partial x}
=(1, 0, \frac{\partial f}{\partial x}, \frac{\partial g}{\partial
x})$ , $e_2=\frac{\partial F}{\partial y} =(0, 1, \frac{\partial
f}{\partial y}, \frac{\partial g}{\partial y})$ and
$v_3=(-\frac{\partial f}{\partial x}, -\frac{\partial f}{\partial
y}, 1, 0)$ and $v_4=(-\frac{\partial g}{\partial x},
-\frac{\partial g}{\partial y}, 0, 1)$. Then $\{e_1, e_2, v_3,
v_4\}$ is a basis of $M$. The metric of $\Sigma$ in this basis is
\begin{eqnarray*}
(g_{ij})=\left(\begin{array}{lr} 1+(\frac{\partial f}{\partial
x})^2+(\frac{\partial g}{\partial x})^2 &\frac{\partial
f}{\partial x}\frac{\partial f}{\partial y}+\frac{\partial
g}{\partial x}\frac{\partial g}{\partial y}\\ \frac{\partial
f}{\partial x}\frac{\partial f}{\partial y}+\frac{\partial
g}{\partial x}\frac{\partial g}{\partial y} &1+(\frac{\partial
f}{\partial y})^2+(\frac{\partial g}{\partial y})^2
\end{array}\right),
\end{eqnarray*} and the inverse matrix is
\begin{eqnarray*}
(g^{ij})=\frac{1}{\det(g)}\left(\begin{array}{lr}
1+(\frac{\partial f}{\partial y})^2+(\frac{\partial g}{\partial
y})^2 &-\frac{\partial f}{\partial x}\frac{\partial f}{\partial
y}-\frac{\partial g}{\partial x}\frac{\partial g}{\partial y}\\
-\frac{\partial f}{\partial x}\frac{\partial f}{\partial
y}-\frac{\partial g}{\partial x}\frac{\partial g}{\partial y}
&1+(\frac{\partial f}{\partial x})^2+(\frac{\partial g}{\partial
x})^2
\end{array}\right).
\end{eqnarray*} It is easy to see that
\begin{eqnarray}\label{cos1}\cos\alpha &=&\frac{\omega(e_1, e_2)}{\sqrt{\det(g)}}
=\frac{\langle Je_1, e_2\rangle}{\sqrt{\det(g)}}\nonumber\\ &=&
\frac{1+\frac{\partial f}{\partial x}\frac{\partial g}{\partial
y}-\frac{\partial g}{\partial x}\frac{\partial f}{\partial
y}}{\det(g)},
\end{eqnarray} and
\begin{eqnarray}\label{sin}
\sin^2\alpha =\frac{\det(g)-\langle Je_1, e_2\rangle^2}{\det(g)}=
\frac{(\frac{\partial f}{\partial x}-\frac{\partial g}{\partial
y})^2+(\frac{\partial g}{\partial x}+\frac{\partial f}{\partial
y})^2}{\det(g)}.
\end{eqnarray}

We express the Euler-Lagrangian equation (\ref{cri}) of $L$
explicitly in the following. We know that,
\begin{eqnarray*}
\nabla\cos\alpha &=& (g^{11}\frac{\partial \cos\alpha}{\partial
x}+g^{12}\frac{\partial\cos\alpha}{\partial
y})e_1+(g^{12}\frac{\partial \cos\alpha}{\partial
x}+g^{22}\frac{\partial\cos\alpha}{\partial y})e_2 \\
&=&\frac{1}{\det(g)}(g_{22}\frac{\partial\cos\alpha}{\partial
x}-g_{12}\frac{\partial\cos\alpha}{\partial
y})e_1+\frac{1}{\det(g)}(g_{11}\frac{\partial\cos\alpha}{\partial
y}-g_{12}\frac{\partial\cos\alpha}{\partial x})e_2 \\ &:=&A e_1+B
e_2.
\end{eqnarray*} Set $c=1+\frac{\partial f}{\partial x}\frac{\partial g}{\partial
y}-\frac{\partial g}{\partial x}\frac{\partial f}{\partial y}$,
$a=\frac{\partial g}{\partial x}+\frac{\partial f}{\partial y}$,
$b=\frac{\partial f}{\partial x}-\frac{\partial g}{\partial y}$.
Then we have,\allowdisplaybreaks
\begin{eqnarray*}
(J(J\nabla\cos\alpha)^\top)^\bot &=& (J(A \frac{\langle Je_1,
e_2\rangle}{|e_2|^2}e_2+B \frac{\langle Je_2,
e_1\rangle}{|e_1|^2}e_1))^\bot \\
&=&(J(\frac{A\cos\alpha\sqrt{\det(g)}}{g_{22}}e_2-\frac{B
\cos\alpha\sqrt{\det{g}}}{g_{11}}e_1))^\bot \\ &=&
\frac{A\cos\alpha\sqrt{\det(g)}}{g_{22}}(\frac{\langle Je_2,
v_3\rangle}{|v_3|^2}
v_3+\frac{\langle Je_2, v_4\rangle}{|v_4|^2} v_4)\\
&&-\frac{B\cos\alpha\sqrt{\det(g)}}{g_{11}}(\frac{\langle Je_1,
v_3\rangle}{|v_3|^2} v_3+\frac{\langle Je_1, v_4\rangle}{|v_4|^2} v_4)\\
&=&\frac{1}{|v_3|^2}
[\frac{A\cos\alpha\sqrt{\det(g)}}{g_{22}}b+\frac{B\cos\alpha\sqrt{\det(g)}}{g_{11}}a]v_3\\
&=&\frac{1}{|v_4|^2}
[\frac{A\cos\alpha\sqrt{\det(g)}}{g_{22}}a-\frac{B\cos\alpha\sqrt{\det(g)}}{g_{11}}b]v_4,
\\ &=&\frac{\cos\alpha}{|v_3|^2\sqrt{\det(g)}}[(\frac{\partial\cos\alpha}{\partial
x}-\frac{g_{12}}{g_{22}}\frac{\partial \cos\alpha}{\partial
y})b\\
&&+(\frac{\partial\cos\alpha}{\partial
y}-\frac{g_{12}}{g_{11}}\frac{\partial\cos\alpha}{\partial
x})a]v_3\\
&&+\frac{\cos\alpha}{|v_4|^2\sqrt{\det(g)}}[(\frac{\partial\cos\alpha}{\partial
x}-\frac{g_{12}}{g_{22}}\frac{\partial \cos\alpha}{\partial
y})a\\
&&-(\frac{\partial\cos\alpha}{\partial
y}-\frac{g_{12}}{g_{11}}\frac{\partial\cos\alpha}{\partial
x})b]v_4.
\end{eqnarray*} On the other hand, we have
\begin{eqnarray*}
\cos^3\alpha
H&=&\frac{\cos^3\alpha}{|v_3|^2}[g^{11}\langle\frac{\partial^2
F}{\partial x^2}, v_3\rangle +g^{12}\langle\frac{\partial^2
F}{\partial x\partial y}, v_3\rangle
+g^{22}\langle\frac{\partial^2 F}{\partial y^2}, v_3]\rangle v_3
\\&&+\frac{\cos^3\alpha}{|v_4|^2}[g^{11}\langle\frac{\partial^2
F}{\partial x^2}, v_3\rangle +g^{12}\langle\frac{\partial^2
F}{\partial x\partial y}, v_3\rangle
+g^{22}\langle\frac{\partial^2 F}{\partial y^2}, v_4]\rangle v_4\\
&=&\frac{\cos^3\alpha}{|v_3|^2\det(g)}(g_{22}\frac{\partial^2
f}{\partial x^2}-2g_{12}\frac{\partial^2 f}{\partial x\partial
y}+g_{11}\frac{\partial^2 f}{\partial y^2})v_3\\
&&+\frac{\cos^3\alpha}{|v_4|^2\det(g)}(g_{22}\frac{\partial^2
g}{\partial x^2}-2g_{12}\frac{\partial^2 g}{\partial x\partial
y}+g_{11}\frac{\partial^2 g}{\partial y^2})v_4.
\end{eqnarray*} Therefore (\ref{cri}) can be written as a systems,
\begin{eqnarray}\label{cri1}
&&\frac{c^2}{(\det(g))^2}(g_{22}\frac{\partial^2 f}{\partial
x^2}-2g_{12}\frac{\partial^2 f}{\partial x\partial
y}+g_{11}\frac{\partial^2 f}{\partial
y^2})\nonumber\\&&-\frac{1}{\sqrt{\det(g)}}[\frac{\partial\cos\alpha}{\partial
x}(b-\frac{g_{12}}{g_{11}}a)+\frac{\partial\cos\alpha}{\partial
y}(a-\frac{g_{12}}{g_{22}}b)]=0\nonumber\\&&\frac{c^2}{(\det(g))^2}(g_{22}\frac{\partial^2
g}{\partial x^2}-2g_{12}\frac{\partial^2 g}{\partial x\partial
y}+g_{11}\frac{\partial^2 g}{\partial
y^2})\nonumber\\&&-\frac{1}{\sqrt{\det(g)}}[\frac{\partial\cos\alpha}{\partial
x}(a+\frac{g_{12}}{g_{11}}b)-\frac{\partial\cos\alpha}{\partial
y}(b+\frac{g_{12}}{g_{22}}a)]=0.
\end{eqnarray} It is not hard to check that,\allowdisplaybreaks
\begin{eqnarray*}
\frac{\partial\cos\alpha}{\partial x}
&=&\frac{1}{(\sqrt{\det(g)})^3}[ \frac{\partial^2 f}{\partial
x^2}(\frac{\partial g}{\partial y}-\frac{\partial f}{\partial
x}+\frac{\partial f}{\partial y}\frac{\partial g}{\partial
x}\frac{\partial g}{\partial y}+\frac{\partial f}{\partial
x}\frac{\partial f}{\partial y}\frac{\partial g}{\partial x}\\
&&+\frac{\partial g}{\partial y}(\frac{\partial g}{\partial
x})^2+\frac{\partial g}{\partial y}(\frac{\partial f}{\partial
y})^2+(\frac{\partial g}{\partial y})^3-\frac{\partial f}{\partial
x}(\frac{\partial g}{\partial y})^2)\\ &+& \frac{\partial^2
g}{\partial x^2}(-\frac{\partial f}{\partial y}-\frac{\partial
g}{\partial x}+\frac{\partial f}{\partial x}\frac{\partial
f}{\partial y}\frac{\partial g}{\partial y}-\frac{\partial
f}{\partial x}\frac{\partial g}{\partial x}\frac{\partial
g}{\partial y}\\ &&-\frac{\partial f}{\partial y}(\frac{\partial
f}{\partial x})^2-\frac{\partial f}{\partial y}(\frac{\partial
g}{\partial y})^2-(\frac{\partial f}{\partial y})^3-\frac{\partial
g}{\partial x}(\frac{\partial f}{\partial
y})^2)\\&+&\frac{\partial^2 f}{\partial x\partial
y}(-\frac{\partial g}{\partial x}-\frac{\partial f}{\partial
y}+\frac{\partial f}{\partial x}\frac{\partial g}{\partial
x}\frac{\partial g}{\partial y}-\frac{\partial f}{\partial
x}\frac{\partial f}{\partial y}\frac{\partial g}{\partial y}\\
&&-\frac{\partial g}{\partial x}(\frac{\partial f}{\partial
x})^2-\frac{\partial g}{\partial x}(\frac{\partial g}{\partial
y})^2-(\frac{\partial g}{\partial x})^3-\frac{\partial f}{\partial
y}(\frac{\partial g}{\partial x})^2)\\ &+&\frac{\partial^2
g}{\partial x\partial y}(\frac{\partial f}{\partial
x}-\frac{\partial g}{\partial y}+\frac{\partial f}{\partial
x}\frac{\partial f}{\partial y}\frac{\partial g}{\partial
x}+\frac{\partial f}{\partial
y}\frac{\partial g}{\partial x}\frac{\partial g}{\partial y}\\
&&+\frac{\partial f}{\partial x}(\frac{\partial g}{\partial
x})^2+\frac{\partial f}{\partial x}(\frac{\partial f}{\partial
y})^2+(\frac{\partial f}{\partial x})^3-\frac{\partial g}{\partial
y}(\frac{\partial f}{\partial x})^2)]\\
&=&\frac{1}{(\sqrt{\det(g)})^3}[ \frac{\partial^2 f}{\partial
x^2}(g_{12}a-g_{22}b)+\frac{\partial^2 g}{\partial
x^2}(-g_{22}a-g_{12}b)\\&&+\frac{\partial^2 f}{\partial x\partial
y}(-g_{11}a+g_{12}b)+\frac{\partial^2 g}{\partial x\partial
y}(g_{11}b+g_{12}a)].
\end{eqnarray*} Similarly, one checks
\begin{eqnarray*}
\frac{\partial\cos\alpha}{\partial y}
&=&\frac{1}{(\sqrt{\det(g)})^3}\frac{\partial^2 f}{\partial
y^2}(-g_{11}a+g_{12}b)+\frac{\partial^2 g}{\partial
y^2}(g_{11}b+g_{12}a)\\ &&+\frac{\partial^2 f}{\partial x\partial
y}(g_{12}a-g_{22}b)+\frac{\partial^2 g}{\partial x\partial
y}(-g_{22}a-g_{12}b).
\end{eqnarray*} Putting these two equations into (\ref{cri1}), then
(\ref{cri1}) can be simplified as,
\begin{eqnarray*}
&&\frac{\partial^2 f}{\partial
x^2}(g_{22}c^2-g_{12}ab+\frac{g_{12}^2}{g_{11}}a^2+g_{22}b^2-\frac{g_{22}g_{12}}{g_{11}}ab)\\&&+
\frac{\partial^2 f}{\partial x\partial
y}(-2g_{12}c^2+g_{11}ab-2g_{12}a^2-2g_{12}b^2+\frac{g_{12}^2}{g_{11}}ab+\frac{g_{12}^2}{g_{22}}ab+g_{22}ab)
\\&&+\frac{\partial^2 f}{\partial
y^2}(g_{11}c^2+g_{11}a^2-g_{12}ab-\frac{g_{11}g_{12}}{g_{22}}ab+\frac{g_{12}^2}{g_{22}}b^2)
\\ &&+\frac{\partial^2 g}{\partial
x^2}(g_{22}ab-\frac{g_{22}g_{12}}{g_{11}}a^2+g_{12}b^2-\frac{g_{12}^2}{g_{11}}ab)
\\ &&+\frac{\partial^2 g}{\partial x\partial
y}(-g_{11}b^2+\frac{g_{12}^2}{g_{11}}a^2+g_{22}a^2-\frac{g_{12}^2}{g_{22}}b^2)
\\ &&+\frac{\partial^2 g}{\partial
y^2}(-g_{11}ab+\frac{g_{11}g_{12}}{g_{22}}b^2-g_{12}a^2+\frac{g_{12}^2}{g_{22}}ab)=0,
\end{eqnarray*} and
\begin{eqnarray*}
&&\frac{\partial^2 g}{\partial
x^2}(g_{22}c^2+g_{22}a^2+\frac{g_{12}^2}{g_{11}}b^2+g_{12}ab+\frac{g_{22}g_{12}}{g_{11}}ab)\\&&+
\frac{\partial^2 g}{\partial x\partial
y}(-2g_{12}c^2-g_{11}ab-2g_{12}a^2-2g_{12}b^2-\frac{g_{12}^2}{g_{11}}ab-\frac{g_{12}^2}{g_{22}}ab-g_{22}ab)
\\&&+\frac{\partial^2 g}{\partial
y^2}(g_{11}c^2+g_{11}b^2+g_{12}ab+\frac{g_{11}g_{12}}{g_{22}}ab+\frac{g_{12}^2}{g_{22}}a^2)
\\ &&+\frac{\partial^2 f}{\partial
x^2}(-g_{12}a^2+\frac{g_{22}g_{12}}{g_{11}}b^2+g_{22}ab-\frac{g_{12}^2}{g_{11}}ab)
\\ &&+\frac{\partial^2 f}{\partial x\partial
y}(g_{11}a^2-\frac{g_{12}^2}{g_{11}}b^2-g_{22}b^2+\frac{g_{12}^2}{g_{22}}a^2)
\\ &&+\frac{\partial^2 f}{\partial
y^2}(-g_{11}ab-\frac{g_{11}g_{12}}{g_{22}}a^2+g_{12}b^2+\frac{g_{12}^2}{g_{22}}ab)=0.
\end{eqnarray*} For simplicity, we write the systems as
\begin{eqnarray*}
A_{11}\frac{\partial^2 f}{\partial x^2}+A_{12}\frac{\partial^2
f}{\partial x\partial y}+A_{22}\frac{\partial^2 f}{\partial
y^2}+B_{11}\frac{\partial^2 g}{\partial
x^2}+B_{12}\frac{\partial^2 g}{\partial x\partial
y}+B_{22}\frac{\partial^2 g}{\partial y^2}
&=&0,\\C_{11}\frac{\partial^2 f}{\partial
x^2}+C_{12}\frac{\partial^2 f}{\partial x\partial
y}+C_{22}\frac{\partial^2 f}{\partial y^2}+D_{11}\frac{\partial^2
g}{\partial x^2}+D_{12}\frac{\partial^2 g}{\partial x\partial
y}+D_{22}\frac{\partial^2 g}{\partial y^2} &=&0,
\end{eqnarray*}
where $A_{ij}$, $B_{ij}$, $C_{ij}$ and $D_{ij}$ are defined
clearly ($i,j=1,2$). So the symbol of the systems is
\begin{eqnarray*}
\sigma:=\left(\begin{array}{lr}
A_{11}\xi^2+A_{12}\xi\eta+A_{22}\eta^2
&B_{11}\xi^2+B_{12}\xi\eta+B_{22}\eta^2\\
C_{11}\xi^2+C_{12}\xi\eta+C_{22}\eta^2
&D_{11}\xi^2+D_{12}\xi\eta+D_{22}\eta^2
\end{array}\right).
\end{eqnarray*} A simple computation yields
\begin{eqnarray*}
\det(\sigma)&=&
\frac{c^4}{g_{22}g_{11}}(g_{11}^3g_{22}\eta^4+g_{22}^3g_{11}\xi^4+2g_{11}^2g_{22}^2\xi^2\eta^2
\\&&+4g_{12}^2g_{11}g_{22}\xi^2 \eta^2-4g_{12}g_{22}g_{11}^2
\xi\eta^3-4g_{12}g_{11}g_{22}^2\xi^3\eta)\\
&&+\frac{c^2(a^2+b^2)}{g_{11}g_{22}}(g_{11}^3g_{22}\eta^4
+g_{11}^2g_{12}^2\eta^4+g_{22}^3g_{11}\xi^4+g_{22}^2g_{12}^2\xi^4
\\&&+2g_{11}^2g_{22}^2\xi^2\eta^2+10g_{12}^2g_{11}g_{22}\xi^2\eta^2-6g_{12}g_{22}g_{11}^2\xi\eta^3 \\ &&
-6g_{12}g_{11}g_{22}^2\xi^3\eta -2g_{12}^3
g_{11}\xi\eta^3-2g_{12}^3g_{22}\xi^3\eta)\\
&=& c^4(g_{22}\xi^2+g_{11}\eta^2-2g_{12}\xi\eta)^2\\
&&+c^2(a^2+b^2)(g_{22}\xi^2+g_{11}\eta^2-3g_{12}\xi\eta)^2\\
&&+\frac{c^2(a^2+b^2)}{g_{11}g_{22}}(g_{11}^2g_{12}^2\eta^4+g_{22}^2g_{12}^2\xi^4+g_{12}^2g_{11}g_{22}\xi^2\eta^2\\
&&~~~~~~-2g_{12}^3g_{11}\xi\eta^3-2g_{12}^3g_{22}\xi^3\eta)
\end{eqnarray*}

Using the inequality $\det(g)=g_{11}g_{22}-g_{12}^2>0$, we have
\begin{eqnarray*}
\det(\sigma)&>& c^4(g_{22}\xi^2+g_{11}\eta^2-2g_{12}\xi\eta)^2\\
&&+c^2(a^2+b^2)(g_{22}\xi^2+g_{11}\eta^2-3g_{12}\xi\eta)^2\\
&&+\frac{c^2(a^2+b^2)g_{12}^2}{g_{11}g_{22}}(g_{11}\eta^2+g_{22}\xi^2-g_{12}\xi\eta)^2
\\ &\geq&0.
\end{eqnarray*}

Therefore the equation (\ref{cri}) is elliptic. This proves the
theorem.

 \hfill Q. E. D.

\section{Equations of the K\"ahler angle of the symplectic critical surfaces}
In the sequel, we always choose the orthonormal basis $\{e_1, e_2,
v_3, v_4\}$ on $M$ along $\Sigma$ such that $\{e_1, e_2\}$ are the
basis of $\Sigma$ and $\omega$ takes the form as (\ref{e2}), and
the complex structure $J$ on $M$ takes the form as (\ref{e3}).

Let $T\Sigma, N\Sigma$ be the tangent bundle and the normal bundle
of $\Sigma$ in $M$ respectively. The second fundamental form $A:
T\Sigma\times T\Sigma\mapsto N\Sigma$ is defined by $A(X,
Y)=(\bar\nabla_X Y)^\bot$ for any tangent vector fields $X, Y$.
The operator $B: T\Sigma\times N\Sigma\mapsto T\Sigma$ is defined
by $B(X, N)=(\bar\nabla_X N)^\top$, $N\in N\Sigma$. Here $()^\top$
denotes the projection from $TM$ onto $T\Sigma$ and $()^\bot$
denote the projection onto $N\Sigma$. Evidently,
$$\langle A(X, Y), N\rangle=-\langle Y, B(X, N)\rangle.
$$

\begin{proposition}\label{pro2}
Let $M$ be a compact K\"ahler surface with K\"ahler form $\omega$
and $J$ be the complex structure compatible with $\omega$ on $M$.
If $\Sigma$ is a closed symplectic surface which is smoothly
immersed in $M$ with the K\"ahler angle $\alpha$, then
\begin{eqnarray}\label{cos2}
\Delta\cos\alpha &=& \cos\alpha(-|h^3_{1k}-h^4_{2k}|^2-
|h^4_{1k}+h^3_{2k}|^2)\nonumber
\\ &&+\sin\alpha(H^4_{,1}+H^3_{,2})-\frac{\sin^2\alpha}{\cos\alpha}
(K_{1212}+K_{1234}).
\end{eqnarray} where $K$ is the curvature operator of $M$ and
$H^\alpha_{,i}=\langle\bar\nabla_{e_i}^N H, v_\alpha\rangle$.
\end{proposition}
{\it Proof.} It is evident that
\begin{eqnarray}\label{e4}
\Delta\cos\alpha &=&\Delta\frac{\omega(e_1,
e_2)}{\sqrt{\det{(g)}}}= \Delta\omega(e_1, e_2)-
\frac{1}{2}\cos\alpha\Delta g_{ij} g^{ij}.
\end{eqnarray} Using the property that $\bar\nabla\omega=0$, we obtain that,
\begin{eqnarray*}
\Delta\omega(e_1, e_2) &=&\bar\nabla_{e_k}\bar\nabla_{e_k}\omega(e_1, e_2) \\
&=&\omega(\bar\nabla_{e_k}\bar\nabla_{e_k} e_1,
e_2)-\omega(\bar\nabla_{e_k}\bar\nabla_{e_k} e_2, e_1)+2\omega(
{\bar\nabla_{e_k} e_1, \bar\nabla_{e_k} e_2}) \\ &=&
\omega(\bar\nabla_{e_k}(\nabla_{e_k} e_1+A(e_k, e_1)), e_2 )-
\omega(\bar\nabla_{e_k}(\nabla_{e_k} e_2+A(e_k, e_2)), e_1 )\\&&+2
\omega(A(e_k, e_1), A(e_k, e_2))\\
&=&\omega(\nabla_{e_k}\nabla_{e_k} e_1, e_2)+\omega(A(e_k,
\nabla_{e_k} e_1), e_2)+\omega(\bar\nabla_{e_k} A(e_k, e_1), e_2)\\
&&-\omega(\nabla_{e_k}\nabla_{e_k} e_2, e_1)+\omega(A(e_k,
\nabla_{e_k} e_2), e_1)-\omega(\bar\nabla_{e_k} A(e_k, e_2), e_1)
\\ &&+2\omega(A(e_k, e_1), A(e_k, e_2))\\
&=&\cos\alpha\langle\nabla_{e_k}\nabla_{e_k} e_1,
e_1\rangle+\cos\alpha\langle\nabla_{e_k}\nabla_{e_k} e_2,
e_2\rangle \\&&+\omega(\bar\nabla_{e_k} A(e_k, e_1),
e_2)-\omega(\bar\nabla_{e_k} A(e_k, e_2), e_1)\\ &&+2\omega(A(e_k,
e_1), A(e_k, e_2)).
\end{eqnarray*} It is not hard to check that
\begin{eqnarray*}
\frac{1}{2}\cos\alpha\Delta g_{ij} g^{ij}&=&\frac{1}{2}\cos\alpha
\Delta\langle e_i, e_j\rangle g^{ij} \\
&=&\cos\alpha\langle\nabla_{e_k}\nabla_{e_k} e_i, e_j\rangle
g^{ij}\\ &=&\cos\alpha\langle\nabla_{e_k}\nabla_{e_k} e_1,
e_1\rangle+\cos\alpha\langle\nabla_{e_k}\nabla_{e_k} e_2,
e_2\rangle.
\end{eqnarray*} Putting the last two identities into (\ref{e4}) and by (\ref{e2}), we
obtain, \allowdisplaybreaks
\begin{eqnarray}\label{e9}
\Delta\cos\alpha &=&\omega(\bar\nabla_{e_k} A(e_k, e_1),
e_2)-\omega(\bar\nabla_{e_k} A(e_k, e_2), e_1)\nonumber\\
&&+2\omega(A(e_k, e_1), A(e_k, e_2))\nonumber \\
&=&\omega(\bar\nabla_{e_k} (h^\alpha_{1k}v_\alpha),
e_2)-\omega(\bar\nabla_{e_k} (h^\alpha_{2k} v_\alpha),
e_1)+2\omega(h^\alpha_{1k}v_\alpha, h^\beta_{2k}v_\beta)\nonumber \\
&=& \omega(h^\alpha_{k1,k}v_\alpha-h^\alpha_{k1}h^\alpha_{kl}e_l,
e_2)-\omega(h^\alpha_{k2,k}v_\alpha-h^\alpha_{k2}h^\alpha_{kl}e_l,
e_1)\nonumber\\&&+2\omega(h^\alpha_{1k}v_\alpha, h^\beta_{2k}v_\beta)\nonumber\\
&=&
\cos\alpha(-(h^\alpha_{1k})^2-(h^\alpha_{2k})^2+2h^3_{1k}h^4_{2k}-2h^4_{1k}h^3_{2k}) \nonumber\\
&&+\omega(h^\alpha_{kk,1}v_\alpha-K_{\alpha k1k}v_\alpha, e_2)\nonumber\\
&&-\omega(h^\alpha_{kk,2}v_\alpha-K_{\alpha k2k}v_\alpha, e_1) \nonumber \\
&=&
\cos\alpha(-(h^\alpha_{1k})^2-(h^\alpha_{2k})^2+2h^3_{1k}h^4_{2k}-2h^4_{1k}h^3_{2k})
\nonumber\\ &&+\sin\alpha(H^4,1+H^3,2)\nonumber\\
&&-\sin\alpha(K_{4k1k}+K_{3k2k}) \nonumber \\&=&
\cos\alpha(-(h^\alpha_{1k})^2-(h^\alpha_{2k})^2+2h^3_{1k}h^4_{2k}-2h^4_{1k}h^3_{2k})
\nonumber\\ &&+\sin\alpha(H^4_{,1}+H^3_{,2})\nonumber\\
&&-\sin\alpha(K_{1213}-K_{1224}) .
\end{eqnarray} Since $J$ is integrable, we have
\begin{eqnarray*}
K_{1212}&=&K(e_1, e_2, Je_1, Je_2)=K(e_1, e_2, \cos\alpha
e_2+\sin\alpha v_3, -\cos\alpha e_1-\sin\alpha v_4) \\
&=& \cos^2\alpha K_{1212}-\sin^2\alpha
K_{1234}+\sin\alpha\cos\alpha (K_{1213}-K_{1224}).
\end{eqnarray*} We therefore obtain
\begin{equation}\label{e7}
\sin\alpha(K_{1212}+K_{1234})=\cos\alpha(K_{1213}-K_{1224}).
\end{equation} Then adding (\ref{e7}) into (\ref{e9}), we get
that,
\begin{eqnarray*}
\Delta\cos\alpha &=&
\cos\alpha(-(h^\alpha_{1k})^2-(h^\alpha_{2k})^2+2h^3_{1k}h^4_{2k}-2h^4_{1k}h^3_{2k})
\\ &&+\sin\alpha(H^4_{,1}+H^3_{,2})-\frac{\sin^2\alpha}{\cos\alpha}(K_{1212}+K_{1234}).
\end{eqnarray*}
This proves the proposition.

\hfill Q. E. D.

\begin{lemma}\label{cor1}
Suppose that $M$ is a K\"ahler surface. We have
\begin{equation}\label{e8}
Ric(Je_1, e_2) =Ric(e_1, Je_2) =\frac{1}{\cos\alpha}
(K_{1212}+K_{1234}).
\end{equation}
\end{lemma}
{\it Proof.} From the Bianchi identity we see that,
\begin{eqnarray*}
Ric(Je_1, e_2)&=& \sum_{A=1}^4K(Je_1, e_A, e_2, e_A) \\
&=& \sum_{A=1}^4K(Je_1, e_A, Je_2, Je_A)\\ &=&-\sum_{A=1}^4K(Je_1,
Je_2, Je_A, e_A)-\sum_{A=1}^4K(Je_1, Je_A, e_A, Je_2) \\
&=&\sum_{A=1}^4K(e_1, e_2, e_A, Je_A)-\sum_{A=1}^4K(Je_1, Je_A, e_2, Je_A)\\
&=&\sum_{A=1}^4K(e_1, e_2, e_A, Je_A)-Ric(Je_1, e_2),
\end{eqnarray*} where we have used the fact that $\{Je_A\}$ are
also orthonormal basis of $M$. Using (\ref{e7}) we get that,
\begin{eqnarray*}
Ric(Je_1, e_2) &=&\frac{1}{2} K(e_1, e_2, e_A, Je_A)\nonumber \\
&=&\sin\alpha(K_{1213}-K_{1224})+\cos\alpha(K_{1212}+K_{1234})\nonumber
\\ &=&\frac{1}{\cos\alpha} (K_{1212}+K_{1234}).
\end{eqnarray*}
Since $Ric(e_1, Je_2) =-Ric(Je_2, e_1)=-\frac{1}{2}K(e_2, e_1,
e_A, Je_A)=\frac{1}{2}K(e_1, e_2, e_A, Je_A),$ it follows that
$$Ric(e_1, Je_2)=Ric(Je_1, e_2).$$

This completes the proof of the lemma.

\hfill Q. E. D.

\begin{theorem}
Suppose that $M$ is K\"ahler surface and $\Sigma$ is a symplectic
critical surface in $M$ with K\"ahler angle $\alpha$, then
$\cos\alpha$ satisfies,
\begin{eqnarray}\label{e12}
\Delta\cos\alpha
&=&\frac{3\sin^2\alpha-2}{\cos\alpha}|\nabla\alpha|^2-\cos\alpha
\sin^2\alpha(K_{1212}+K_{1234})\nonumber\\ &=&
\frac{3\sin^2\alpha-2}{\cos\alpha}|\nabla\alpha|^2-\cos^2\alpha
\sin^2\alpha Ric(Je_1, e_2) .
\end{eqnarray}
\end{theorem}

{\it Proof.} If $\Sigma$ is a symplectic critical surface, then
$H=\frac{\sin^2\alpha}{\cos^2\alpha}V$. It is easy to check that,
\begin{eqnarray*}
(h^3_{1k}-h^4_{2k})^2+(h^4_{1k}+h^3_{2k})^2 &=&
|H|^2+2|V|^2+2H\cdot V \\ &=&
(\frac{\sin^4\alpha}{\cos^4\alpha}+2+2\frac{\sin^2\alpha}{\cos^2\alpha})|V|^2
\\ &=& \frac{1+\cos^4\alpha}{\cos^4\alpha} |\nabla\alpha|^2,
\end{eqnarray*} and
\begin{eqnarray*}
H^4_{, 1}+H^3_{,2}
&=&\partial_1(\frac{\sin^2\alpha}{\cos^2\alpha}\partial_1\alpha)
+\partial_2(\frac{\sin^2\alpha}{\cos^2\alpha}\partial_2\alpha) \\
&=&
\frac{\sin^2\alpha}{\cos^2\alpha}\Delta\alpha+2\frac{\sin\alpha}{\cos^3\alpha}
|\nabla\alpha|^2 \\
&=&\frac{\sin\alpha}{\cos^2\alpha}(-\Delta\cos\alpha-\cos\alpha|\nabla\alpha|^2)
+2\frac{\sin\alpha}{\cos^3\alpha} |\nabla\alpha|^2.
\end{eqnarray*} Putting these two equations into (\ref{cos2}), we
obtain that,
\begin{eqnarray*}
\Delta\cos\alpha
&=&\frac{-1-\cos^4\alpha+2\sin^2\alpha-\sin^2\alpha\cos^2\alpha}{\cos^3\alpha}
|\nabla\alpha|^2-\frac{\sin^2\alpha}{\cos^2\alpha}\Delta\cos\alpha
\\ &&-\frac{\sin^2\alpha}{\cos\alpha}(K_{1212}+K_{1234}).
\end{eqnarray*} Lemma \ref{cor1} implies that,
\begin{eqnarray*}
\Delta\cos\alpha &=& \frac{3\sin^2\alpha-2}{\cos\alpha}
|\nabla\alpha|^2-\cos\alpha\sin^2\alpha(K_{1212}+K_{1234}) \\ &=&
\frac{3\sin^2\alpha-2}{\cos\alpha}
|\nabla\alpha|^2-\cos^2\alpha\sin^2\alpha Ric(Je_1, e_2).
\end{eqnarray*}

This completes the proof of the theorem.

\hfill Q. E. D.

\begin{corollary}\label{cor2}
If $M$ is a K\"ahler-Einstein surface with scalar curvature $R$,
and $\Sigma$ is a symplectic critical surface in $M$ with K\"ahler
angle $\alpha$, we have
\begin{eqnarray*}
\Delta\cos\alpha &=& \frac{3\sin^2\alpha-2}{\cos\alpha}
|\nabla\alpha|^2-R\cos^3\alpha\sin^2\alpha ,
\end{eqnarray*}
if in addition, we assume that $R\geq 0$, then
 the symplectic critical surface $\Sigma$ is a holomorphic curve.
\end{corollary}

{\it Proof.} Suppose that $M$ is a K\"ahler-Einstein surface with
scalar curvature $R$, we have
\begin{eqnarray*}
Ric(Je_1, e_2)=\cos\alpha R_{22}=\cos\alpha R,
\end{eqnarray*}
The identity in the corollary follows. The second statement of the
corollary follows from the maximum principle. This proves the
corollary.

 \hfill Q. E. D.

\section{Topology of the symplectic critical surfaces}
In this section we will analyze the topology properties of the
symplectic critical surfaces. At a point $p\in \Sigma$ with
$\alpha(p)=0$ the tangent plane $T_p\Sigma$ of $M$ at $p$ is a
complex line in $T_{F(p)}M$. So such a point is called a complex
tangent point.  We recall some equations obtained by Wolfson in
\cite {W2} (also see \cite {W1}). We write the metric of $M$ as
$$
ds^2=\sum_{\beta=-1,1}\omega_\beta\bar{\omega}_\beta,
$$
the induced metric on $\Sigma$ can be written as
$$
ds_\Sigma^2=\phi\circ\bar\phi,
$$
where $\phi$ is a complex valued $1$-form defined up to a complex
factor of norm one, furthermore, one can assume that
$$
\omega_{1}=\cos(\frac{\alpha}{2})\phi,~~\omega_{-1}=\sin(\frac{\alpha}{2})\bar\phi
,
$$
where $\alpha$ is the K\"ahler angle.

Suppose that the complex second fundamental form of $\Sigma$ in
$M$ is
$$
II^C=a\phi^2+2b\phi\bar\phi +c\bar\phi^2.
$$

Relative to the coframe field $\omega_{-1}$, $\omega_1$, there is
a unitary connection $\omega_{\beta\gamma}$ which satisfies
$$
d\omega_\beta=\omega_{\beta\gamma}\wedge\omega_\gamma,~~\omega_{\beta\gamma}+\omega_{\gamma\beta}=0.
$$

We set
$$
\cos(\alpha/2)\omega_1+\sin(\alpha/2)\bar\omega_{-1}=\theta_1+\sqrt{-1}\theta_2,
$$
$$
\sin(\alpha/2)\omega_1-\cos(\alpha/2)\bar\omega_{-1}=\theta_3+\sqrt{-1}\theta_4,
$$
where $\theta_k$, $k=1,\cdots,4$, is an orthonormal coframe of the
Riemannian structure of $M$. So, along $\Sigma$ we have
$$
\sin(\alpha/2)\omega_1-\cos(\alpha/2)\bar\omega_{-1}=0,
$$
it follows that (\cite {W1} (1.6), \cite {W2} (2.18))
\begin{equation}\label{wolfson01}
\frac{1}{2}(d\alpha +\sin\alpha
(\omega_{-1\bar{-1}}+\omega_{1\bar1}))=a\phi+b\bar\phi .
\end{equation}

The relation between the real second fundamental form and the
complex second fundamental form is given in \cite {W2} (Section
2),
\begin{eqnarray*}
&&\left(\begin{array}{lr} 1  &1\\ \sqrt{-1}
&-\sqrt{-1}\end{array}\right)\left(\begin{array}{lr} a  &b\\
b &c\end{array}\right)\left(\begin{array}{lr} 1 &\sqrt{-1}\\
1 &-\sqrt{-1}\end{array}\right)\\&&=\left(\begin{array}{lr} h^3_{11}  &h^3_{12}\\
h^3_{12}
&h^3_{22}\end{array}\right)+\sqrt{-1}\left(\begin{array}{lr} h^4_{11}  &h^4_{12}\\
h^4_{12} &h^4_{22}\end{array}\right).
\end{eqnarray*}
We therefore have
\begin{equation}\label{wolfson02}
b=\frac{1}{4}(H^3+\sqrt{-1}H^4).
\end{equation}
Using the equation (\ref{wolfson01}), (\ref{wolfson02}) and
(\ref{cri}), we see that, on the symplectic critical surface we
have,
$$
\frac{\partial \sin\alpha}{\partial \bar{\zeta}}=(\sin\alpha )h,
$$
where $h$ is a smooth complex function, $\zeta$ is a local complex
coordinate on $\Sigma$. By Bers' \cite {B} result, we have

\begin{proposition}\label{pro4}
A non holomorphic symplectic critical surface in a K\"ahler
surface has at most finite complex points.
\end{proposition}
Set $$ g(\alpha)=\ln(\sin^2\alpha).$$ Then using the equation
(\ref{e12}), we obtain
\begin{eqnarray}\label{e13}
\Delta g(\alpha)
&=&-2|\nabla\alpha|^2-2\frac{\cos\alpha}{\sin^2\alpha}\Delta\cos\alpha
-4\frac{\cos^2\alpha}{\sin^2\alpha}|\nabla\alpha|^2 \nonumber\\
&=& -4|\nabla\alpha|^2+2\cos^2\alpha(K_{1212}+K_{1234}).
\end{eqnarray} This equation is valid away from the complex
tangent points of $M$. By the Gauss equation and Ricci equation,
we have,
\begin{eqnarray*}
R_{1212}&=&K_{1212}+h^\alpha_{11}h^\alpha_{22}-(h^\alpha_{12})^2\\
R_{1234}&=&K_{1234}+h^3_{1k}h^4_{2k}-h^4_{1k}h^3_{2k},
\end{eqnarray*} where $R_{1212}$ is the curvature of $T\Sigma$ and $R_{1234}$ is the
curvature of $N\Sigma$. Adding these two equations together, we
get that,
\begin{eqnarray*}
K_{1212}+K_{1234} &=& R_{1212}+R_{1234}-\frac{1}{2}|H|^2 \\
&&+\frac{1}{2}((h^3_{1k}-h^4_{2k})^2+(h^4_{1k}+h^3_{2k})^2)\\ &=&
R_{1212}+R_{1234}+|V|^2+H\cdot V \\ &=&
R_{1212}+R_{1234}+\frac{1}{\cos^2\alpha}|\nabla \alpha|^2.
\end{eqnarray*} Thus,
\begin{eqnarray*}
R_{1212}+R_{1234}&=&\frac{1}{2}\frac{\Delta
g(\alpha)}{\cos^2\alpha}+\frac{|\nabla\alpha|^2}{\cos^2\alpha}.
\end{eqnarray*} Integrating the above equality over $\Sigma$ we
have,
\begin{eqnarray*}
2\pi(\chi(T\Sigma)+\chi(N\Sigma))&=&-2\pi
P-\int_{\Sigma}\frac{|\nabla\alpha|^2}{\cos^2\alpha} d\mu_\Sigma,
\end{eqnarray*} where $\chi(T\Sigma)$ is the Euler characteristic
of $\Sigma$ and $\chi(N\Sigma)$ is the Euler characteristic of the
normal bundle of $\Sigma$ in $M$, $P$ is the sum of the orders of
complex tangent points. We therefore proved the following theorem.
\begin{theorem}\label{top01}
Let $\Sigma$ be a non holomorphic symplectic critical surface in a
K\"ahler surface $M$. Let $P$ denote the sum of the orders of
complex tangent points. Then
$$\chi(T\Sigma)+\chi(N\Sigma)=-
P-\frac{1}{2\pi}\int_{\Sigma}\frac{|\nabla\alpha|^2}{\cos^2\alpha}
d\mu_\Sigma.
$$
\end{theorem}

\begin{remark}
Let $\Sigma$ be a symplectic critical surface in a K\"ahler
surface $M$. Then
$$\chi(T\Sigma)+\chi(N\Sigma)\leq 0.
$$
\end{remark}

Similarly, we can show,

\begin{theorem}\label{top02}
Let $\Sigma$ be a non holomorphic symplectic critical surface in a
K\"ahler surface $M$, then
$$F^\ast c_1(M)[\Sigma]=- P-\frac{1}{2\pi}\int_{\Sigma}\frac{|\nabla\alpha|^2}{\cos^3\alpha}d\mu_\Sigma.
$$
where $c_1(M)$ is the first Chern class of $M$ and [$\Sigma$] is
the fundamental homology class of $\Sigma$.
\end{theorem}

{\it Proof.} By (\ref{e13}) we also get that,
\begin{eqnarray*}
\Delta g(\alpha)&=& -4|\nabla\alpha|^2+2\cos^3\alpha Ric(Je_1,
e_2).
\end{eqnarray*} Note that $Ric(Je_1, e_2) d\mu_\Sigma$ is the pulled back Ricci 2-form of $M$ by the
immersion $F$ to $\Sigma$, i.e,
$$ F^\ast (Ric^M)=Ric(Je_1, e_2) d\mu_\Sigma.$$ Thus,
$$F^\ast (Ric^M)=(\frac{1}{2}\frac{\Delta g(\alpha)}{\cos^3\alpha}
+2\frac{|\nabla\alpha|^2}{\cos^3\alpha})d\mu_\Sigma.
$$ Integrating it over $\Sigma$, we obtain that,
$$2\pi F^\ast c_1(M)[\Sigma]=-2\pi P
-\int_{\Sigma}\frac{|\nabla\alpha|^2}{\cos^3\alpha}d\mu_\Sigma.$$
This proves the theorem. \hfill Q. E. D.
\begin{remark}
If $\Sigma$ is a symplectic critical surface in a K\"ahler surface
$M$., then
$$F^\ast c_1(M)[N]\leq 0.$$
\end{remark}

\begin{corollary}\label{cor3}
Suppose that $\Sigma$ is a symplectic critical surface in a
K\"ahler surface $M$. Then
$$\chi(T\Sigma)+\chi(N\Sigma)\geq  F^\ast c_1(M)[N].
$$ and the equality holds if and only if $\Sigma$ is a holomorphic curve.
\end{corollary}

{\it Proof.} It is easy to see that $F^\ast
c_1(M)[N]\leq\chi(T\Sigma)+\chi(N\Sigma)$ and the equality holds
on the holomorphic curves. We need only prove the necessary part.
Suppose that $F^\ast c_1(M)[N]=\chi(T\Sigma)+\chi(N\Sigma)$, then
$$\int_{\Sigma}\frac{|\nabla\alpha|^2}{\cos^2\alpha}d\mu_\Sigma
=\int_{\Sigma}\frac{|\nabla\alpha|^2}{\cos^3\alpha}d\mu_\Sigma.
$$ That means,
$$ \int_\Sigma \frac{|\nabla\alpha|^2}{\cos^2\alpha} (
\frac{1}{\cos\alpha}-1)d\mu_\Sigma=0,
$$ which implies that $\cos\alpha=1$, i.e, $\Sigma$ is a
holomorphic curve.

 \hfill Q. E. D.

Recall that the degree of a map $F: \Sigma\to M$ is defined by
$$ \deg\Sigma=\frac{1}{\pi}\int_{\Sigma} F^\ast\omega.
$$

\begin{corollary}\label{cor4}
Suppose that $\Sigma$ is a non holomorphic symplectic critical
surface in a K\"ahler-Einstein surface $M$ with scalar curvature
$-R>0$. Then
$$
\frac{1}{2}R\deg \Sigma
=P+\frac{1}{2\pi}\int_{\Sigma}\frac{|\nabla\alpha|^2}{\cos^3\alpha}d\mu.
$$
\end{corollary}

{\it Proof.} It simply follows from (\ref{e13}). The detail is
left to readers.

\hfill Q. E. D.

Let $g$ be the genus of $\Sigma$, $I_\Sigma$ the self-intersection
number of $\Sigma$, $D_\Sigma$ the number of double points of
$\Sigma$. Then
$$
\chi(\Sigma)=2-2g,
$$
$$
\chi(\nu)=I_\Sigma-2D_\Sigma.
$$

Set
$$
C_1(\Sigma)=C_1(M)([\Sigma]),
$$
we have
\begin{theorem}\label{th3}
Suppose that $\Sigma$ is a symplectic critical surface in a
K\"ahler surface $M$. Then
$$
2-2g-
c_1(\Sigma)+I_\Sigma-2D_\Sigma=\frac{1}{2\pi}\int_{\Sigma}\frac{(1-\cos\alpha)|\nabla\alpha|^2}{\cos^3\alpha}d\mu.
$$
\end{theorem}
{\it Proof.} It simply follows from Theorem \ref{top01} and
Theorem \ref{top02}. The detail is left to readers.

\hfill Q. E. D.

\section{The gradient flow}

In this section we consider the gradient flow of the function $L$,
i. e,
\begin{equation}\label{main01}
\frac{dF}{dt}=\cos^2\alpha
H-\frac{1}{\cos\alpha}(J(J\nabla\cos\alpha)^\top)^\bot.
\end{equation}
We set
$$
\vec{f}=\cos^2\alpha
H-\frac{1}{\cos\alpha}(J(J\nabla\cos\alpha)^\top)^\bot.
$$
It is clear that, if $\vec{f}=0$, $\Sigma$ is a symplectic
critical surface.

By the first variational formula of the functional $L$ (Theorem
\ref{th1}), we see that, along the flow,
\begin{eqnarray}\label{dec}
\frac{dL}{dt}&=&-2\int_{\Sigma}\frac{1}{\cos^3\alpha}|\cos^2\alpha
H-\frac{1}{\cos\alpha}(J(J\nabla\cos\alpha)^\top)^\bot|^2d\mu\nonumber\\
&=&-2\int_{\Sigma}\frac{1}{\cos^3\alpha}|\vec{f}|^2d\mu.
\end{eqnarray}

By Theorem \ref{th2}, we know that the equation (\ref{main01}) is
a parabolic equation, and the short time existence can be shown by
a standard argument. We set $\Sigma_t=F(\Sigma, t)$ with
$\Sigma_0=\Sigma$.

Using the same local frame as that in Section 2, the equation
(\ref{main01}) can be written as
\begin{equation}\label{main1}
\frac{dF}{dt}=\cos^2\alpha\vec{H}-\sin^2\alpha\vec{V}\equiv\vec{f}.
\end{equation}

We compute the evolution of the area element of $\Sigma_t$ along
the flow.
\begin{lemma}\label{l1}
\begin{equation}\label{e10}\frac{d}{dt}d\mu_t=(-\cos^2\alpha
|H|^2+\sin^2\alpha \vec V\cdot\vec H)d\mu_t. \end{equation}
\end{lemma}

{\it Proof.} It is easy to check that,
\begin{eqnarray*}
\frac{\partial}{\partial t}g_{ij} &=&\frac{\partial}{\partial
t}\langle\frac{\partial F}{\partial x^i},\frac{\partial
F}{\partial x^j}\rangle=2\langle\frac{\partial \vec{f}}{\partial
x^i}, \frac{\partial F}{\partial x^j}\rangle \\
&=&-2\langle\vec{f}, \frac{\partial^2 F}{\partial x^i\partial
x^j}\rangle \\ &=& 2(-\cos^2\alpha H^\alpha+\sin^2\alpha V^\alpha)
h^\alpha_{ij}.
\end{eqnarray*} Here we have used the fact that $\nabla_{e_i} e_j=0$ at the
fixed point. Therefore, we have \begin{eqnarray*}
\frac{d}{dt}d\mu_t &=& -\vec{f}\cdot\vec H d\mu_t \\
&=&(-\cos^2\alpha |H|^2+\sin^2\alpha \vec V\cdot\vec H) d\mu_t.
\end{eqnarray*}

\hfill Q. E. D.

Now we derive the evolution equation of $\cos\alpha$ along the
flow (\ref{main1}), which can be seen as a start point of the
study of the flow.
\begin{theorem}\label{pro3} Let $M$ be a K\"ahler surface. Assume
that $\alpha$ is the K\"ahler angle of $\Sigma_t$ which evolves by
the flow (\ref{main1}). Then $\cos\alpha$ satisfies the equation
\begin{eqnarray}\label{cos}
(\frac{d}{dt}-\Delta)\cos\alpha &=&\cos^3\alpha
(|h^3_{1k}-h^4_{2k}|^2+|h^3_{2k}+h^4_{1k}|^2)
+\cos^2\alpha\sin^2\alpha Ric(Je_1,
e_2)\nonumber\\&&+\cos\alpha\sin^2\alpha
|H|^2-\cos\alpha\sin^2\alpha|V+H|^2,
\end{eqnarray} where $\{e_1, e_2, v_3, v_4\}$ is an orthonormal
basis of $T_pM$ such that $\omega, J$ take the form (\ref{e2}),
(\ref{e3}).
\end{theorem}

{\it Proof.} Using the fact that $\bar\nabla\omega=0$ and by Lemma
\ref{l1}, we have,
\begin{eqnarray*}
\frac{\partial}{\partial t}\cos\alpha &=&\frac{\partial}{\partial
t}\frac{\omega(e_1,
e_2)}{\sqrt{\det(g_t)}}=\omega(\bar\nabla_{e_1}\vec{f},
e_2)-\omega(\bar\nabla_{e_2}\vec{f},
e_1)-\frac{1}{2}\cos\alpha\frac{\partial}{\partial t}g_{ij} g^{ij}
\\ &=&\omega(\bar\nabla_{e_1}\vec{f},
e_2)-\omega(\bar\nabla_{e_2}\vec{f}, e_1)+\cos\alpha\vec{f}\cdot
H.
\end{eqnarray*} By breaking $\bar\nabla_{e_1}\vec{f}$ and $ \bar\nabla_{e_2}\vec{f}$
into normal and tangent parts, we get that,
\begin{eqnarray*}
\omega(\bar\nabla_{e_1}\vec{f},
e_2)-\omega(\bar\nabla_{e_2}\vec{f}, e_1)
&=&\omega(\bar\nabla_{e_1}^N \vec{f},
e_2)-\omega(\bar\nabla_{e_2}^N \vec{f}, e_1)\\
&&+\omega(\bar\nabla_{e_1}^T \vec{f},
e_2)-\omega(\bar\nabla_{e_2}^N \vec{f}, e_1)\\
&=&\omega(\bar\nabla_{e_1}^N \vec{f},
e_2)-\omega(\bar\nabla_{e_2}^N \vec{f},
e_1)\\ &&+\omega(B(e_1, \vec{f}), e_2)-\omega(B(e_2, \vec{f}), e_1)\\
&=&\omega(\bar\nabla_{e_1}^N \vec{f},
e_2)-\omega(\bar\nabla_{e_2}^N \vec{f}, e_1)\\
&&+\cos\alpha(\langle
B(e_1, \vec{f}), e_1\rangle+\langle B(e_2, \vec{f}), e_2\rangle)\\
&=&\omega(\bar\nabla_{e_1}^N \vec{f},
e_2)-\omega(\bar\nabla_{e_2}^N \vec{f},
e_1)-\cos\alpha\vec{f}\cdot\vec{ H}.
\end{eqnarray*} Combining these two identities we obtain that,
\begin{eqnarray}\label{e11}
\frac{\partial}{\partial t}\cos\alpha &=&
\omega(\bar\nabla_{e_1}^N \vec{f}, e_2)-\omega(\bar\nabla_{e_2}^N
\vec{f}, e_1)\nonumber\\ &=&\omega(\bar\nabla_{e_1}^N(\cos^2\alpha
H-\sin^2\alpha V), e_2)-\omega(\bar\nabla_{e_2}^N(\cos^2\alpha
H-\sin^2\alpha V), e_1)\nonumber\\
&=&\cos^2\alpha(\omega(\bar\nabla_{e_1}^N H,
e_2)-\omega(\bar\nabla_{e_2}^N H, e_1))\nonumber\\
&&+\omega(\partial_1\cos^2\alpha H,
e_2)-\omega(\partial_2\cos^2\alpha H, e_1)\nonumber \\
&&-\omega(\partial_1\sin^2\alpha V,
e_2)+\omega(\partial_2\sin^2\alpha V, e_1)\nonumber \\
&&-\sin^2\alpha(\omega(\bar\nabla_{e_1}^N V,
e_2)-\omega(\bar\nabla_{e_2}^N V, e_1))\nonumber\\
&:=&
I+II+III+IV.
\end{eqnarray}
By Theorem \ref{pro2} we have
\begin{eqnarray*}
I &=& \cos^2\alpha\sin\alpha(H^4_{,1}-H^3_{,2}) \\ &=&\cos^2\alpha
(\Delta\cos\alpha+\cos\alpha |h^3_{1k}-h^4_{2k}|^2+\cos\alpha
|h^3_{2k}+h^4_{1k}|^2+\sin^2\alpha Ric(Je_1, e_2)).
\end{eqnarray*} It is clear that
\begin{eqnarray*}
II &=&
H^4\sin\alpha\partial_1\cos^2\alpha+H^3\sin\alpha\partial_2\cos^2\alpha
\\ &=& -2\sin^2\alpha\cos\alpha(H^4\partial_1\alpha+H^3\partial_2\alpha) \\ &=&
-2\cos\alpha\sin^2\alpha \vec H\cdot\vec V.
\end{eqnarray*} From the definition of $V$, we can see that,
\begin{eqnarray*}
III &=& -\sin\alpha\partial_1\sin^2\alpha\partial_1\alpha-
\sin\alpha\partial_2\sin^2\alpha\partial_2\alpha \\
&=&-2\sin^2\alpha\cos\alpha|\nabla\alpha|^2 \\
&=&-2\sin^2\alpha\cos\alpha |V|^2 .
\end{eqnarray*}
Similarly, one obtains
\begin{eqnarray*}
IV &=&-\sin^2\alpha\omega(\bar\nabla_{e_1}^N (\partial_2\alpha
v_3+\partial_1\alpha v_4),
e_2)+\sin^2\alpha\omega(\bar\nabla_{e_2}^N (\partial_2\alpha
v_3+\partial_1\alpha v_4), e_1)\\
&=&-\sin^2\alpha(\omega(\partial_1\partial_1\alpha v_4,
e_2)-\omega(\partial_2\partial_2\alpha v_3, e_1))\\
&=&-\sin^3\alpha\Delta\alpha\\
&=&\sin^2\alpha\Delta\cos\alpha+\sin^2\alpha\cos\alpha
|\nabla\alpha|^2\\ &=&
\sin^2\alpha\Delta\cos\alpha+\sin^2\alpha\cos\alpha |V|^2.
\end{eqnarray*} Putting these equations into (\ref{e11}), we
can obtain that,
\begin{eqnarray*}
\frac{\partial}{\partial t}\cos\alpha &=&
\Delta\cos\alpha+\cos^3\alpha |h^3_{1k}-h^4_{2k}|^2+\cos^3\alpha
|h^3_{2k}+h^4_{1k}|^2\\ &&+\cos^2\alpha\sin^2\alpha Ric(Je_1,
e_2)-2\cos\alpha\sin^2\alpha \vec H\cdot\vec
V-\cos\alpha\sin^2\alpha |V|^2.
\end{eqnarray*} This proves the theorem.

\hfill Q. E. D.

\begin{theorem}\label{pro3ke} Let $M$ be a K\"ahler-Einstein
surface with scalar curvature $R$.
Assume that $\alpha$ is the K\"ahler angle of $\Sigma_t$ which
evolves by the flow (\ref{main1}). Then $\cos\alpha$ satisfies the
equation
\begin{eqnarray}\label{coske}
(\frac{d}{dt}-\Delta)\cos\alpha &=&\cos^3\alpha
(|h^3_{1k}-h^4_{2k}|^2+|h^3_{2k}+h^4_{1k}|^2)
+R\cos^3\alpha\sin^2\alpha \nonumber\\
&&+\cos\alpha\sin^2\alpha
|H|^2-\cos\alpha\sin^2\alpha|V+H|^2,
\end{eqnarray} where $\{e_1, e_2, v_3, v_4\}$ is an orthonormal
basis of $T_pM$ such that $\omega, J$ take the form (\ref{e2}),
(\ref{e3}). And consequently, if $\Sigma$ is symplectic, along the
flow (\ref{main1}), at each time $t$,  $\Sigma_t$ is symplectic.
\end{theorem}

\end{document}